\documentstyle[prb,aps]{revtex}
\begin{document}
\title{Expressing the Field Operator in Terms of Currents and Densities}
\author{Girish S. Setlur}
\address{Department of Physics and Materials Research Laboratory,\\
 University of Illinois at Urbana-Champaign , Urbana Il 61801}
\maketitle
    
\maketitle

\begin{abstract}
 A formula is written down for the annhilation operator(bose or 
 fermi) in terms of the corresponding observable bilinears namely
 currents and densities. The Fock space representation of these
 formulas is clarified. A conjecture is written down that 
 shows that such a correspondence is unique. 
\end{abstract}
\section{Introduction}
  In the late 70's and early 80's, attempts were made to write down field
 theories that describe scalar mesons in terms of observables like 
 currents and densities rather than the creation and annhilation operators.  
 The motivation for doing this stems from the fact a theory cast directly
 in terms of observables was more physically intuitive than the more traditional
 approach based on raising and lowering operators on the Fock space.
 This attempt however, raised a number of technical questions, among them
 was how to make sense of the various identities connecting say the kinetic
 energy density to the currents and particle densities and so on. 
 Elaborate mathematical machinery was erected by the authors who started
 this line of research\cite{Sharp} to address these issues. However, it seems
 gaps still remain especially with regard to the crucial question of how
 one goes about writing down a formula for the annhilation operator (fermi or
 bose) alone in terms of bilinears like currents and densities. The 
 bilinears in question namely currents and densities satisfy a closed algebra
 known as the current algebra \cite{Sharp}. This algebra is insensitive
 to the nature of the statistics of the underlying fields. On the other hand,
 if one desires information about single-particle properties, it is necessary
 to relate the annhilation operator (whose commutation rules determine the
 statistics) to bilinears like currents and densities. That such a
 correspondence is possible was demonstrated by Goldin, Menikoff and Sharp
 \cite{Sharp}. However their methods are somewhat abstract and difficult to
 follow especially for someone trained exclusively in Physics.  This 
 is the main motivation for my article which is aimed primarily at Physicists.
 This allows me the luxury of sacrificing some mathematical rigor and its place
 make plausibility arguments in the hope of drawing in a wider audience.
  The attempts made here are partly based on the work of Goldin et.al.
  \cite{Sharp}  Ligouri and
  Mintchev on generalised statistics\cite{Lig} and the series by
  Reed and Simon\cite{Reed}. This article is organised as follows. 
 In the next couple of paragraphs, we prove a lemma that is going to be 
 used repeatedly in the following sections. In the next section, we
 write the DPVA ansatz\cite{Setlur3}(this was introduced in an earlier 
 preprint and it stands for "Density Phase Variable Approach")
 in the Fock space language. This also means relating
 the canonical conjugate of the density to observables like currents
 and densities. In the end we write down the main conjecture which 
 when proven will vindicate the DPVA ansatz\cite{Setlur3}.
 A partial proof(no rigor) of the
 conjucture is given elsewhere\cite{setlurexact}. Now for the Lemma.
\newline
{\bf{Lemma}}
 Let $ {\mathcal{F}} $ be smooth a function from a bounded
 subset of the real line on to the set of reals. Also let
 $ f $ and $ g $ be smooth functions from some bounded subset of
 $ {\mathcal{R}}^{d} $ to reals.
 Let us further assume that
 the range of these functions are such that it is always possible to
 find compositions such as $ {\mathcal{F}}\mbox{ }o\mbox{ }f $
 and they will also be smooth functions with sufficiently big domains.
 The claim is this, since they possess fourier transforms (they are
 well-behaved) If,
\begin{equation}
{\mathcal{F}}(\mbox{  }f({\bf{x}})\mbox{  }) = g({\bf{x}})
\end{equation}
and,
\begin{equation}
f({\bf{x}}) = \sum_{ {\bf{q}} }{\tilde{f}}_{ {\bf{q}} }
\mbox{  }e^{ i\mbox{ }{\bf{q.x}} }
\end{equation}
\begin{equation}
g({\bf{x}}) = \sum_{ {\bf{k}} }{\tilde{g}}_{ {\bf{k}} }
\mbox{  }e^{ i\mbox{ }{\bf{k.x}} }
\end{equation}
then the following also holds,
\begin{equation}
[ {\mathcal{F}}(\mbox{ }\sum_{ {\bf{q}} }
{\tilde{f}}_{ {\bf{q}} }T_{ -{\bf{q}} }({\bf{k}}) \mbox{ })]
\mbox{ }\delta_{ {\bf{k}}, {\bf{0}} } = {\tilde{g}}_{ {\bf{k}} }
\end{equation}
where $ T_{ {\bf{q}} }({\bf{k}}) = exp({\bf{q}}.\nabla_{ {\bf{k}} }) $.
Here the operator $ T_{ {\bf{q}} }({\bf{k}}) $ acts on the $ {\bf{k}} $
 in the Kronecker delta on the extreme right, and every time it translates
 the $ {\bf{k}} $ by an amount $ {\bf{q}} $.
\newline
{\bf{Proof}}
\newline
Proof is by brute force expansion.
We know,
\begin{equation}
{\mathcal{F}}(y) = \sum_{n=0}^{\infty} \frac{ {\mathcal{F}}^{(n)}(0) }{n!}
\mbox{ }y^{n}
\end{equation}
therefore,
\[
{\mathcal{F}}( \mbox{  }f({\bf{x}}) \mbox{  })
 = {\mathcal{F}}(0) + \sum_{n=1}^{\infty} \frac{ {\mathcal{F}}^{(n)}(0) }{n!}
\mbox{ }\sum_{ \{ {\bf{q}}_{i} \} }{\tilde{f}}_{ {\bf{q}}_{1} }
{\tilde{f}}_{ {\bf{q}}_{2} }...{\tilde{f}}_{ {\bf{q}}_{n} }
exp(i(\sum_{i=1}^{n}{\bf{q}}_{i}).{\bf{x}})
\]
\begin{equation}
= \sum_{ {\bf{k}} }e^{i\mbox{ }{\bf{k.x}} }{\tilde{g}}_{ {\bf{k}} }
\end{equation}
This means (take inverse fourier transform),
\[
 {\mathcal{F}}(0)\delta_{ {\bf{k}}, 0 }
 + \sum_{n=1}^{\infty} \frac{ {\mathcal{F}}^{(n)}(0) }{n!}
\mbox{ }\sum_{ \{ {\bf{q}}_{i} \} }{\tilde{f}}_{ {\bf{q}}_{1} }
{\tilde{f}}_{ {\bf{q}}_{2} }...{\tilde{f}}_{ {\bf{q}}_{n} }
\delta_{ ( {\bf{k}} - \sum_{i=1}^{n}{\bf{q}}_{i} ), {\bf{0}} }
\]
\begin{equation}
= {\tilde{g}}_{ {\bf{k}} }
\end{equation}
This may also be cleverly rewritten as
\[
 {\mathcal{F}}(0)\delta_{ {\bf{k}}, 0 }
 + \sum_{n=1}^{\infty} \frac{ {\mathcal{F}}^{(n)}(0) }{n!}
\mbox{ }\sum_{ \{ {\bf{q}}_{i} \} }{\tilde{f}}_{ {\bf{q}}_{1} }
{\tilde{f}}_{ {\bf{q}}_{2} }...{\tilde{f}}_{ {\bf{q}}_{n} }
T_{ -{\bf{q}}_{1} }({\bf{k}})T_{ -{\bf{q}}_{2} }({\bf{k}})
...T_{ -{\bf{q}}_{n} }({\bf{k}})
\delta_{  {\bf{k}}, {\bf{0}} }
\]
\begin{equation}
= {\tilde{g}}_{ {\bf{k}} }
\end{equation}
and therefore,
\begin{equation}
{\tilde{g}}_{ {\bf{k}} } = [{\mathcal{F}}
(\mbox{  }\sum_{ {\bf{q}} }{\tilde{f}}_{ {\bf{q}} }
T_{ -{\bf{q}} }({\bf{k}}) \mbox{  } )]\mbox{  }\delta_{ {\bf{k}}, {\bf{0}} }
\end{equation}
and the {\bf{Proof}} is now complete.
                        
\section{Fock Space Representation}
 In earlier preprints\cite{Setlur3},
 we introduced the so-called DPVA ansatz(see Appendix also
 Ref[\onlinecite{Setlur3}])
 that related the annhilation operator to the canonical conjugate of
 the density distribution. In this section, we try to formulate this
 correspondence in the Fock space language.            
 We start off with some preliminaries. Let $ {\mathcal{H}} $ be an
 infinite dimensional separable Hilbert Space. We know from
 textbooks that such a space possesses a countable orthonormal
 basis $ {\mathcal{B}} = \{ w_{i}; i \in {\mathcal{Z}} \} $.
 Here, $ {\mathcal{Z}} $ is the set of all integers.
 Thus$ {\mathcal{H}} $ = Set of all linear combinations of vectors chosen
 from $ {\mathcal{B}} $. We construct the tensor product of two such
 spaces
\[
 {\mathcal{H}}^{ {\small{\bigotimes}} 2}
 = {\mathcal{H}} {\bigotimes} {\mathcal{H}}
\]
 This is defined to be the dual space of the space of all bilinear
 forms on the direct sum. In plain English this means something like this.
 Let $ f \in {\mathcal{H}} $ and $ g \in {\mathcal{H}} $ define the
 object $ f {\small{\bigotimes}} g $
 to be that object which acts as shown below.
 Let $ < v, w > $ be an element of the Cartesian product
 $ {\mathcal{H}} \times {\mathcal{H}} $.
\[
 f {\small{\bigotimes}} g < v, w > = (f, v) (g, w)
\]
Here, $ (f, v) $ stands for the inner product of $ f $ and $ v $.
Define also the inner product of two $ f {\small{\bigotimes}} g $ and
  $ f^{'} {\small{\bigotimes}} g^{'} $
\[
(f {\small{\bigotimes}} g,  f^{'} {\small{\bigotimes}}
 g^{'}) = (f,  f^{'} ) (g, g^{'})
\]
 Construct the space of all finite linear combinations of
 objects such as $ f {\small{\bigotimes}} g $ with different choices for
 $ f $ and $ g $. Lump them all into a set. You get a vector space.
 It is still not the vector space
 $ {\mathcal{H}}^{{\small{\bigotimes 2}}} $.
 Because the space of all finite linear combinations of
 objects such as $ f {\small{\bigotimes}} g $ is not complete. Not every
 Cauchy sequence converges. Complete the space by appending the limit points
 of all Cauchy sequences from the space of all finite linear combinations
 of vectors of the type $ f {\small{\bigotimes}} g $.
 This complete space is the
 Hilbert space $ {\mathcal{H}}^{{\small{\bigotimes 2}}} $. Similarly one can
 construct $  {\mathcal{H}}^{{\small{\bigotimes}} n} $ for n = 0, 1, 2, 3, ...
 Where we have set $ {\mathcal{H}}^{0}  = {\mathcal{C}} $
 the set of complex numbers. Define the Fock Space over
 $ {\mathcal{H}} $ as
\[
{\mathcal{F}}({\mathcal{H}}) = {\bigoplus}_{n=0}^{\infty}
 {\mathcal{H}}^{{\small{\bigotimes}} n}
\]
 Physically, each element of it is an ordered collection of wavefunctions
 with different number of particles
\[
 (\varphi_{0}, \varphi_{1}(x_{1}), \varphi_{2}(x_{1}, x_{2}),
 ...,\varphi_{n}, ... )
\]
 is a typical element of $ {\mathcal{F}}({\mathcal{H}}) $. This is the
 Hilbert Space which we shall be working with. Let $ {\mathcal{D}}^{n} $
 be the space of all decomposable vectors.
\[
{\mathcal{D}}^{n} = \{ f_{1} {\small{\bigotimes}}...{\small{\bigotimes}} f_{n}
; f_{i} \in {\mathcal{H}} \}
\]
For each $ f \in {\mathcal{H}} $ define
\[
b(f): {\mathcal{D}}^{n} \rightarrow {\mathcal{D}}^{n-1}, n \geq 1
\]
\[
b^{*}(f): {\mathcal{D}}^{n} \rightarrow {\mathcal{D}}^{n+1}, n \geq 0
\]
defined by
\[
b(f) \mbox{ } f_{1} {\small{\bigotimes}} ...  {\small{\bigotimes}} f_{n}
 = \sqrt{n} (f, f_{1}) f_{2} {\small{\bigotimes}} ...{\small{\bigotimes}}
 f_{n}
\]

\[
b^{*}(f) \mbox{ } f_{1} {\small{\bigotimes}} ...{\small{\bigotimes}} f_{n}
 = \sqrt{n+1}
 f \bigotimes  f_{1} {\small{\bigotimes}} f_{2}
{\small{\bigotimes}} ...  {\small{\bigotimes}} f_{n}
\]
 We also define $ b(f) {\mathcal{H}}^{0}  = 0 $. By linearity we can extend
 the definitions to the space of all finite linear combinations
 of elements of $ {\mathcal{D}}^{n} $ namely
 $ {\mathcal{L}} ({\mathcal{D}}^{n}) $.
 For any $ \varphi \in {\mathcal{L}}({\mathcal{D}}^{n}) $ and
$ \psi \in {\mathcal{L}}({\mathcal{D}}^{n+1}) $
\[
\parallel b(f) \varphi \parallel \leq \sqrt{n} \parallel f \parallel
 \parallel \varphi \parallel
\]
\[
\parallel b^{*}(f) \varphi \parallel \leq \sqrt{n+1} \parallel f \parallel
 \parallel \varphi \parallel
\]
\[
(\psi,  b^{*}(f)\varphi) = ( b(f)\psi, \varphi)
\]
So long as $ \parallel f \parallel < \infty $, $  b(f) $ and $  b^{*}(f)  $
 are bounded operators. An operator $ {\mathcal{O}} $ is said to be bounded
 if
\[
 sup_{\parallel \varphi \parallel = 1}
\mbox{ }
 \parallel {\mathcal{O}} {\mathcal{\varphi}} \parallel \mbox{ }
 < \mbox{ } \infty
\]
 $ {\mathcal{O}} $ is unbounded otherwise. The norm of a bounded operator is
 defined as
\[
\parallel {\mathcal{O}} \parallel =  sup_{\parallel \varphi \parallel = 1}
\mbox{ }
 \parallel {\mathcal{O}} {\mathcal{\varphi}} \parallel \mbox{ }
\]
In order to describe fermions, it is necessary to
 construct orthogonal projectors on $ {\mathcal{F}}({\mathcal{H}}) $.
 In what follows $ c(f) $ will denote a fermi annhilation operator.
 $ c^{*}(f) $ will denote a fermi creation operator. Physically, and naively
 speaking, these are the fermi operators in "momentum space" $ c_{\bf{k}} $
 and $ c^{*}_{\bf{k}} $.
 First define $ P_{-} $ to be the projection operator that projects
 out only the antisymmetric parts of many body wavefunctions.
 For example,
\[
P_{-} f_{1} {\small{\bigotimes}} f_{2} = \frac{1}{2}
 (f_{1} {\small{\bigotimes}} f_{2} - f_{2} {\small{\bigotimes}} f_{1})
\]
We now have
\[
c(f) = P_{-}b(f)P_{-}
\]
\[
c^{*}(f) = P_{-}b^{*}(f)P_{-}
\]
Let us take a more complicated example. Let us find out how
 $ c^{*}(f) c(g) $ acts on a vector
 $ v =  f_{1} {\small{\bigotimes}} f_{2} $.
\[
c^{*}(f) c(g) = P_{-}b^{*}(f)P_{-}P_{-}b(g)P_{-}
\]
\[
c^{*}(f) c(g) = P_{-}b^{*}(f)P_{-}b(g)P_{-}
\]
\[
c^{*}(f) c(g) v = P_{-}b^{*}(f)P_{-}b(g)P_{-} v
\]
\[
P_{-} v = \frac{1}{2 !}( f_{1} {\small{\bigotimes}} f_{2}
 - f_{2} {\small{\bigotimes}} f_{1} )
\]

\[
b(g)P_{-} v = \frac{1}{2 !}{\sqrt{2}}( (g, f_{1}) f_{2} -
 (g, f_{2}) f_{1} )
\]
\[
P_{-}b(g)P_{-} v = \frac{1}{2 !}{\sqrt{2}}( (g, f_{1}) f_{2} -
 (g, f_{2}) f_{1} )
\]
\[
b^{*}(f)P_{-}b(g)P_{-} v = \frac{1}{2 !}{\sqrt{2}}^{2}
( (g, f_{1}) f {\small{\bigotimes}} f_{2} -
 (g, f_{2}) f {\small{\bigotimes}} f_{1} )
\]
\[
c^{*}(f) \mbox{ } c(g) v = (\frac{1}{2 !})^{2}{\sqrt{2}}^{2}
( (g, f_{1})[ f {\small{\bigotimes}}f_{2} - f_{2} {\small{\bigotimes}} f]
 -
 (g, f_{2}) [ f {\small{\bigotimes}} f_{1} - f_{1} {\small{\bigotimes}} f ] )
\]
 Having had a feel for how the fermi operators behave, we are now
 equipped to pose some more pertinent questions.
 Choose a basis
\[
{\mathcal{B}} = \{ w_{i}; i \in {\mathcal{Z}} \}
\]
                
\subsection{Definition of the Fermi Density Distribution}

 Here we would like to capture the notion of the fermi density operator.
 Physicists call it $ \rho(x) = \psi^{*}(x) \psi(x) $. Multiplication
 of two fermi fields at the same point is a tricky business and we would
 like to make more sense out of it. For this we have to set our single
 particle Hilbert Space:
\[
 {\mathcal{H}} = L_{p}({\mathcal{R}}^{3}) {\bigotimes} {\mathcal{W}}
\]
Here, $ L_{p}({\mathcal{R}}^{3}) $ is the space of all periodic functions
 with period $ L $ in each space direction. That is if
 $ u \in  L_{p}({\mathcal{R}}^{3}) $ then
\[
u(x_{1} + L, x_{2}, x_{3}) = u(x_{1}, x_{2}, x_{3})
\]
\[
u(x_{1}, x_{2} + L, x_{3}) = u(x_{1}, x_{2}, x_{3})
\]
\[
u(x_{1}, x_{2}, x_{3} + L) = u(x_{1}, x_{2}, x_{3})
\]

$ {\mathcal{W}} $ is the spin space spanned by two vectors.
 An orthonormal basis for $ {\mathcal{W}} $
\[
\{ \xi_{\uparrow}, \xi_{\downarrow} \}
\]
A typical element of $ {\mathcal{H}} $ is given by
 $ f({\bf{x}}) \bigotimes \xi_{\downarrow} $. A basis for $ {\mathcal{H}} $
 is given by
\[
{\mathcal{B}} =
 \{ \sqrt{ \frac{1}{L^{3}} } exp(i{\bf{q_{n}.x}}) \bigotimes
 \xi_{s}
  ;
 {\bf{n}} = (n_{1}, n_{2}, n_{3}) \in {\mathcal{Z}}^{3},
  s \in \{ \uparrow, \downarrow \} ;
\]
\[
 { \bf{q_{n}} } =
 (\frac{2 \pi n_{1}}{L}, \frac{2 \pi n_{2}}{L}, \frac{2 \pi n_{3}}{L})
 \}
\]
We move on to the definition of the fermi-density distribution.
 The Hilbert Space $ {\mathcal{H}}^{\bigotimes n} $ is the space of
 all n-particle wavefunctions with no symmetry restrictions.
 From this we may construct orthogonal subspaces
\[
{\mathcal{H}}_{+}^{\bigotimes n} = P_{+} {\mathcal{H}}^{\bigotimes n}
\]
\[
{\mathcal{H}}_{-}^{\bigotimes n} = P_{-} {\mathcal{H}}^{\bigotimes n}
\]
 Tensors from $ {\mathcal{H}}_{+}^{\bigotimes n} $ are orthogonal
 to tensors from $ {\mathcal{H}}_{-}^{\bigotimes n} $. The only exceptions are
 when $ n = 0 $ or $ n = 1 $.
\[
{\mathcal{H}}_{+}^{\bigotimes 0} = {\mathcal{H}}_{+}^{\bigotimes 0} =
 {\mathcal{C}}
\]
\[
{\mathcal{H}}_{+}^{\bigotimes 1} = {\mathcal{H}}_{+}^{\bigotimes 1} =
 {\mathcal{H}}
\]
The space $ {\mathcal{H}}_{+}^{\bigotimes n} $ is the space of
 bosonic-wavefunctions and the space $ {\mathcal{H}}_{-}^{\bigotimes n} $
 is the space of fermionic wavefunctions. The definition of the fermi
 density distribution proceeds as follows. Let $ v $ be written as
\[
 v = \sum_{\sigma \in \{ \uparrow, \downarrow \} }a(\sigma)\xi_{\sigma}
\]
 The Fermi density distibution is an operator on the Fock Space,
 given a vector $ f \bigotimes v \in \mathcal{H} $ in the single particle
 Hilbert Space, and a tensor $ \varphi $ in the n-particle subspace of
 of $ \mathcal{F}(\mathcal{H}) $, there exists a corresponding operator
 $ \rho(f \bigotimes v) $ that acts as follows:
\[
[\rho(f \bigotimes v) \varphi ]_{n}
({\bf{x_{1}}}\sigma_{1}, {\bf{x_{2}}}\sigma_{2},
 .... , {\bf{x_{n}}}\sigma_{n})
  = 0
\]
  if $ \varphi \in  {\mathcal{H}}_{+}^{\bigotimes n} $ and
\[
[\rho(f \bigotimes v) \varphi]_{n}
 ({\bf{x_{1}}}\sigma_{1}, {\bf{x_{2}}}\sigma_{2},
 .... , {\bf{x_{n}}}\sigma_{n})
 = \sum_{i=1}^{n} f({\bf{x_{i}}}) a(\sigma_{i})
 \varphi_{n}({\bf{x_{1}}}\sigma_{1}, {\bf{x_{2}}}\sigma_{2},
 .... , {\bf{x_{n}}}\sigma_{n})
\]
when $ \varphi \in  {\mathcal{H}}_{-}^{\bigotimes n} $.
The physical meaning of this abstract operator will become clear
 in the next subsection.

\subsection{Definition of the Canonical Conjugate of the Fermi Density}

 We introduce some notation. Let
$ g = exp(i{ \bf{k_{m}.x} } ) \bigotimes \xi_{r} $ (the square root
 of the volume is not needed in the definition of $ g $ since 
 we want all the operators in the fourier space such as 
$ \psi({ { \bf{k_{m}} } r}) $ to be dimensionless.
\[
\psi({ { \bf{k_{m}} } r}) = c(g)
\]
\[
\rho({ { \bf{k_{m}} } r}) = \rho(g)
\]
This $ \rho({ { \bf{k_{m}} } r})  $ is nothing but the density operator
 in momentum space, familiar to Physicists
\[
 \rho({{\bf{k_{m}}}r}) =
 \sum_{ {\bf{q_{n}}} }c^{\dagger}_{ {\bf{q_{n}+k_{m}}}r }c_{ {\bf{q_{n}}}r }
\]
 We want to define the canonical conjugate of the density operator
 as an operator that maps the Fock space(or a subset thereof)
 on to itself.
 If $ \varphi \in {\mathcal{H}}^{0} = {\mathcal{C}} $ then
\[
X_{ { \bf{q_{m}} } } \varphi = 0
\]
 Let $ \varphi \in  {\mathcal{H}}_{+}^{n} $, $ n = 2, 3, ... $ then
\[
X_{ { \bf{q_{m}} } } \varphi = 0
\]
 The important cases are when
 $ \varphi \in  {\mathcal{H}}_{-}^{n} $, $ n = 2, 3, ... $ or
 if $ \varphi \in  {\mathcal{H}} $. In such a case, we set
 $ n = N_{s}^{0} \neq 0 $ for $ s \in \{ \downarrow, \uparrow \} $.
 Let us introduce some more notation.
 $ N_{s} = \rho_{ { \bf{q = 0} } s } $ is the number operator to be
 distinguished from the c-number $ N_{s}^{0} $. The eigenvalue of
 $ N_{s} $ is $ N_{s}^{0} $ when it acts on a state such as
 $ \varphi \in {\mathcal{H}}_{-}^{n} $. Some more notation.
\[
  \delta \mbox{ } \psi({ { \bf{k_{m}} } s })
  =  \psi({ { \bf{k_{m}} } s }) - { \sqrt{N_{s}^{0}} }
\delta_{ { \bf{k_{m}, 0} } }
\]
 and
\[
 \delta \mbox{ } \rho({ { \bf{k_{m}} } s })
 = \rho({ { \bf{k_{m}} } s }) - N_{s}^{0} \delta_{ { \bf{k_{m}, 0} } }
\]
The Canonical Conjugate of the density distribution in real space
 denoted by $ \Pi_{s}({\bf{x}}) $ is defined as
 follows.
\[
\Pi_{s}(x_{1}+L, x_{2}, x_{3}) = \Pi_{s}(x_{1}, x_{2}, x_{3})
\]
\[
\Pi_{s}(x_{1}, x_{2}+L, x_{3}) = \Pi_{s}(x_{1}, x_{2}, x_{3})
\]
\[
\Pi_{s}(x_{1}, x_{2}, x_{3}+L) = \Pi_{s}(x_{1}, x_{2}, x_{3})
\]
The definition is as follows.
\[
\Pi_{s}({\bf{x}}) = \sum_{ {\bf{q_{m}}}}
 exp(i{\bf{q_{m}.x}})X_{ {\bf{q_{m}}}s }
\]

\[
 X_{ {\bf{q_{m}}}s } =
 i \mbox{ } {\large{ln}} [ ( 1 +
\frac{1}{ { \sqrt{N_{s}^{0}} } }
 \sum_{ { \bf{k_{n}} } } \delta \psi({ \bf{k_{n}} } s)
 T_{- {\bf{k_{n}}} }({\bf{q_{m}}}))( 1 + \sum_{ { \bf{k_{n}} } }
\]
                                                       
\begin{equation}
\frac{1}{ N_{s}^{0} } \delta \rho({ \bf{k_{n}} } s)
 T_{{\bf{k_{n}}} }({\bf{q_{m}}}) )^{-\frac{1}{2}}
exp(-i\sum_{ { \bf{k_{n}} } }\phi([\rho]; {\bf{k_{n}}} s)
T_{{\bf{k_{n}}} }({\bf{q_{m}}})  ) ] \delta_{ {\bf{q_{m}}}, 0}
\label{result}
\end{equation}
\[
\Phi([\rho]; {\bf{x}}s) = \sum_{ {\bf{k_{n}}} }
\phi([\rho];{\bf{k_{n}}} s) exp(-i{\bf{k_{n}.x}})
\]
\[
T_{{\bf{k_{n}}} }({\bf{q_{m}}}) = exp({\bf{k_{n}.\nabla_{q_{m}}}})
\]
 The translation operator translates the $ {\bf{q_{m}}} $
 in the Kronecker delta that appears in the extreme right by
 $ {\bf{k_{n}}} $
 and $ \Phi([\rho]; {\bf{x}}s) $ satisfies a recursion explained in detail
 in the previous manuscript. The logarithm is to be interpreted as
 an expansion around the leading term  which is either $ N_{s}^{0} $ or
 $ { \sqrt{N_{s}^{0}} } $. The question of existence of $ X_{ {\bf{q_{m}}}s } $
 now reduces to demonstrating that this operator (possibly unbounded)
 maps its domain of definition(densely defined in Fock space)
  on to the Fock space. Defining the limit of the series expansion
 is likely to be the major bottleneck in demonstrating the existence
 of $ X_{ {\bf{q_{m}}}s } $. That this is the canonical conjugate of the
 density operator is not at all obvious from the above definition. A
 rigorous proof of that is also likely to be difficult. Considering
 that we arrived at this formula by first postulating the existence of
 $ \Pi_{s}({\bf{x}}) $, it is probably safe to just say
 " it is clear that "  $ \Pi_{s}({\bf{x}}) $, is in fact the canonical
 conjugate of $ \rho $. The way in which the above formula can be
 deduced may be motivated as follows:
\[
 \psi({\bf{x}}\sigma) =
 \frac{ 1 }{ V^{\frac{1}{2}} }
\sum_{ {\bf{k}} }exp(i{\bf{k.x}})\psi({\bf{k}}\sigma)
\]
\[
 = exp(-i \sum_{{\bf{q}}} exp(i{\bf{q.x}}) X_{{\bf{q}}\sigma} )
exp(i \sum_{{\bf{q}}} exp(-i{\bf{q.x}}) \phi( [\rho]; {\bf{q}}\sigma ) )
\]                                         
\begin{equation}
(n_{\sigma} + \frac{1}{V} \sum_{{\bf{q}}\neq 0}
 \rho_{{\bf{q}}\sigma} exp(-i{\bf{q.x}}) )^{\frac{1}{2}}
\label{FF}
\end{equation}
In the above Eq.(~\ref{FF}) ONLY on the right side make the replacements,
\begin{equation}
exp(i{\bf{q.x}}) \rightarrow T_{{\bf{-q}}}({\bf{k}})
\end{equation}
\begin{equation}
exp(-i{\bf{q.x}}) \rightarrow T_{{\bf{q}}}({\bf{k}})
\end{equation}
where
\begin{equation}
T_{{\bf{q}}}({\bf{k}}) = exp({\bf{q}}.\nabla_{{\bf{k}}})
\end{equation}
and append a $ \delta_{ {\bf{k, 0}} } $ on the extreme right.
 Also on the LEFT side of Eqn.(~\ref{FF})
 make the replacement $ \psi({\bf{x}}\sigma) $ by
 $ \psi({\bf{k}}\sigma) $. And you get a formula for
 $  \psi({\bf{k}}\sigma) $.
 In order to get a formula for $ X_{ {\bf{q_{m}} } s } $
 we have to invert the relation and obtain,
\[
\Pi({\bf{x}}s) =
i \mbox{ } ln[  ( \sqrt{N_{s}^{0}} +
\sum_{ {\bf{k_{n}}}}
exp(i{\bf{k_{n}.x}})\delta\psi({\bf{k_{n}}}s) )
\]
\begin{equation}
exp(-i \sum_{{\bf{k_{n}}}} exp(-i{\bf{k_{n}.x}})
 \phi( [\rho]; {\bf{k_{n}}}s ) )(N_{s}^{0} + \sum_{ {\bf{k_{n}}} }
\delta\rho_{{\bf{k_{n}}}s} exp(-i{\bf{k_{n}.x}}) )^{-\frac{1}{2}} ]
\label{canon}
\end{equation}
Make the replacements on the right side of Eqn.(~\ref{canon})
\begin{equation}
exp(i{\bf{k_{n}.x}}) \rightarrow T_{{\bf{-k_{n}}}}({\bf{q_{m}}})
\end{equation}
\begin{equation}
exp(-i{\bf{k_{n}.x}}) \rightarrow T_{{\bf{k_{n}}}}({\bf{q_{m}}})
\end{equation}
where,
and append a $ \delta_{ {\bf{q_{m}, 0}} } $ on the extreme right.
Also on the LEFT side of Eqn.(~\ref{canon}) replace
$ \Pi({\bf{x}}s) $ by $ X_{ {\bf{q_{m}}}s } $. This results in formula
 given in Eqn. (~\ref{result}).

\section{The Canonical Conjugate in Terms of the Current}

We would like to express the field operator in terms of currents and densities.
This is inspired directly by the work of Sharp, Menikoff and Goldin.
To this end let us make the following statements, 
If one insists on having a self-adjoint
 canonical conjugate of the density then one must sacrifice positivity of
 the density. On the other hand if one insists on having a positive density,
 one must sacrifice self-adjointness of the canonical conjugate\cite{Goldpriv}.
 This is
 true even when the underlying single-particle hilbert space is made separable.
 Let us assume that we have decided one way or the other. Then, we
 would like to make a statement that gives us a rigorous way of deciding
 whether it is possible to write down a formula for the field operator in terms
 of currents and densities. To this end, let us do the following,
 first define the current operator rigorously. To Physicists, it is,
\begin{equation}
{\bf{J}}({\bf{x}}) = (\frac{1}{2i})[ \psi^{\dagger}(\nabla\psi)
 - (\nabla\psi)^{\dagger}\psi ]
\end{equation}
 To the Math-community it is an operator similar to the density
\cite{Mintpriv}, given a typical element
 $ f \bigotimes v $ associated with the underlying single-particle
 hilbert space, there is an operator denoted
 by $ J_{s}(f \bigotimes v) $, ( $ s = 1,2,3 $ )
 that acts on a typical tensor from the 
 n-particle subspace of the full Fock space as follows,
\begin{equation}
 [J_{s}(f \bigotimes v)\varphi]_{n}
({\bf{x}}_{1}\sigma_{1}, {\bf{x}}_{2}\sigma_{2}, ... , {\bf{x}}_{n}\sigma_{n})
 = 0
\end{equation}
[6~if $ \varphi \in {\mathcal{H}}^{\bigotimes n}_{+} $,and
and,
\[
 [J_{s}(f \bigotimes v)\varphi]_{n}
({\bf{x}}_{1}\sigma_{1}, {\bf{x}}_{2}\sigma_{2}, ... , {\bf{x}}_{n}\sigma_{n})
\]
\begin{equation}
 = -i\sum_{k = 1}^{n}
\{ f({\bf{x}}_{k})a(\sigma_{k})\nabla^{k}_{s} + \frac{1}{2}
[\nabla^{k}_{s}f({\bf{x}}_{k})]a(\sigma_{k}) \} \varphi_{n}
({\bf{x}}_{1}\sigma_{1}, {\bf{x}}_{2}\sigma_{2}, ... , {\bf{x}}_{n}\sigma_{n})
\end{equation}
if $ \varphi \in {\mathcal{H}}^{\bigotimes n}_{-} $.
For the bosonic current it is the other way around. 
Having done all this, we would now like to write the DPVA ansatz more
 rigorously.
 Now for some notation. As before, let
 $ g = \mbox{ }exp(i{\mbox{ }}{\bf{k}}_{m}.{\bf{x}})\bigotimes \xi_{r} $
(the square root of the volume is not needed as we want all operators 
in momentum space to be dimensionless).
Then as before,
\begin{equation}
\psi({\bf{k}}_{m}r) = c(g)
\end{equation}
\begin{equation}
\rho({\bf{k}}_{m}r) = \rho(g)
\end{equation}
\begin{equation}
\delta \rho({\bf{k}}_{m}r) = \rho({\bf{k}}_{m}r) - N^{0}_{r}
\delta_{ {\bf{k}}_{m}, {\bf{0}} }
\end{equation}
\begin{equation}
j_{s}({\bf{k}}_{m}r) = {\bf{J}}_{s}(g)
\end{equation}
\begin{equation}
\delta j_{s}({\bf{k}}_{m}r) = j_{s}({\bf{k}}_{m}r)
\end{equation}
Having done this, we would like to write down another formula for the canonical
conjugate. Recall that\cite{Setlur3},
\begin{equation}
\nabla\Pi({\bf{x}}\sigma) = (-1/\rho({\bf{x}}\sigma)){\bf{J}}({\bf{x}}\sigma)
+ \nabla\Phi([\rho];{\bf{x}}\sigma) - [-i\mbox{ }\Phi,\nabla\Pi]
\end{equation}
Then we have(bear in mind here that we have distinguished between the c-number
 $ N^{0}_{r} $ and the operator $ \rho({\bf{0}}r) $ whose expectation value
 is $ N^{0}_{r} $).
\begin{equation}
(i\mbox{ }{\bf{q}}_{m})\mbox{  }X_{ {\bf{q}}_{m}r }
 = -(\frac{1}{ N^{0}_{r} })\frac{1}{1 + \frac{1}{N^{0}_{r}}
\sum_{ {\bf{k}}_{n} }
\delta\rho({\bf{k}}_{n}r)T_{ {\bf{k}}_{n} }({\bf{q}}_{m} )}
[ \mbox{  }
\sum_{ {\bf{p}}_{n} }\delta {\bf{j}}({\bf{p}}_{n}r)
T_{ {\bf{p}}_{n} }({\bf{q}}_{m} )
\mbox{   } ] \mbox{   } \delta_{ {\bf{q}}_{m}, {\bf{0}} }
+ \mbox{    }{\bf{F}}([\rho];{\bf{q}}_{m}r)
\end{equation}
where,
\begin{equation}
\sum_{ {\bf{q}}_{m} }exp(i{\mbox{  }}{\bf{q}}_{m}.{\bf{x}})
{\bf{F}}([\rho];{\bf{q}}_{m}r)
 = \nabla\Phi - [-i\mbox{ }\Phi,\nabla\Pi]
\end{equation}
Without loss of generality, we may set $ X_{ {\bf{0}}r } = 0 $, as
 this contributes just a constant phase. For
 $ {\bf{q}}_{m} \neq {\bf{0}} $ 
\[
X_{ {\bf{q}}_{m}r }
 = (\frac{1}{q_{m}^{2}})(\frac{ i }{N^{0}_{r}})
\frac{1}{1 + \frac{1}{N^{0}_{r}}
\sum_{ {\bf{k}}_{n} }
\delta\rho({\bf{k}}_{n}r)T_{ {\bf{k}}_{n} }({\bf{q}}_{m} )}
[ \mbox{  }
\sum_{ {\bf{p}}_{n} }
{\bf{q}}_{m}.\delta {\bf{j}}({\bf{p}}_{n}r)
T_{ {\bf{p}}_{n} }({\bf{q}}_{m} )
\mbox{   } ] \mbox{   } \delta_{ {\bf{q}}_{m}, {\bf{0}} }
\]
\begin{equation}
- \frac{ i\mbox{ }{\bf{q}}_{m}.{\bf{F}}([\rho];{\bf{q}}_{m}r) }
{ q_{m}^{2} }
\end{equation}
Now define an operator which is defined to be the formal expansion
 that the formula itself suggests,
\begin{equation}
{\tilde{\psi}}({\bf{k}}_{n}r) = 
exp(-i{\mbox{  }}
\sum_{ {\bf{q}}_{m}  }T_{ -{\bf{q}}_{m} }({\bf{k}}_{n})
X_{ {\bf{q}}_{m}r })
\mbox{ }exp(i{\mbox{  }}
\sum_{ {\bf{q}}_{m} }T_{ {\bf{q}}_{m} }({\bf{k}}_{n})
\phi([\rho];{\bf{q}}_{m}r)){\mbox{  }}
(N^{0}_{r} + \sum_{ {\bf{q}}_{m} }\delta\rho({\bf{q}}_{m}r)
T_{ {\bf{q}}_{m} }({\bf{k}}_{n}))^{\frac{1}{2}}
\delta_{ {\bf{k}}_{n}, {\bf{0}} }
\end{equation}

\section{The Main Conjecture}

In this section, we would like to write to write down a statement that 
would require a proof. This conjecture when proven will vindicate the
 DPVA ansatz. Let us focus on fermions, the bose case is analogous.
\newline
{\bf{Conjecture}}
\newline
There exists a unique functional $ \Phi([\rho];{\bf{x}}r) $ and a unique odd
(for fermions, even for bosons)
 integer $ m $ such that the following recursion holds, 
\[
\Phi([\{\rho({\bf{y_{1}}}\sigma_{1})
 - \delta({\bf{y_{1}}}-{\bf{x}}^{'})\delta_{\sigma_{1},\sigma^{'}} \} ]
;{\bf{x}}\sigma)
\]
\[
+ \Phi([\rho];{\bf{x^{'}}}\sigma^{'}) - \Phi([\rho];{\bf{x}}\sigma)
\]
\begin{equation}
-\Phi([\{\rho({\bf{y_{1}}}\sigma_{1})
 - \delta({\bf{y_{1}}}-{\bf{x}})\delta_{\sigma_{1},\sigma} \} ]
;{\bf{x^{'}}}\sigma^{'})
 = m\pi
\label{recur}
\end{equation}
 and has the following additional effects.
The domain of definition of $ {\tilde{\psi}}({\bf{k}}_{n}r) $ 
( in which the series expansion converges ), is the same as that of
 $ \psi({\bf{k}}_{n}r) $ and it acts the same way too. In other words,
\begin{equation}
{\tilde{\psi}}({\bf{k}}_{n}r) = \psi({\bf{k}}_{n}r) 
\end{equation}
We know how the ingredients of $ {\tilde{\psi}}({\bf{k}}_{n}r) $ namely
 the current $ {\bf{j}}({\bf{k}}_{n}r) $ and the density
 $ \delta\rho({\bf{q}}_{n}r) $  act on typical elements of the Fock space,
 and we know how $  \psi({\bf{k}}_{n}r) $ acts on the Fock space,
 we just have to show that the complicated $ {\tilde{\psi}}({\bf{k}}_{n}r) $
acts the same way as the simple $ \psi({\bf{k}}_{n}r) $. Moreover, 
 this is true for a unique phase functional $ \Phi $.
\section{Appendix}
 In this section we try to recapture in a compact way the results of
 earlier works\cite{Setlur3},\cite{setlurexact}.
 In particular, the DPVA(which stands for the density
 phase variable approach) ansatz is brought out in a way that is 
 readily appreciated by Physicists. The DPVA ansatz relates the annhilation
 operator to the canonical conjugate of the relevant density distribution.
 It may be written as follows,
\begin{equation}
\psi({\bf{x}}) = exp(-i\mbox{ }\Pi({\bf{x}}) )
exp(i\mbox{ }\Phi([\rho];{\bf{x}}) ) (\rho({\bf{x}}))^{\frac{1}{2}}
\end{equation}
 here $ \Phi $ satisfies the recursion explained earlier. $ \Pi $ 
 is canonically conjugate to the density $ \rho = \psi^{*}\psi $.
 It satisfies canonical equal-time commutation rules with the density.
\begin{equation}
[\Pi({\bf{x}}), \rho({\bf{y}})] = i\mbox{ }\delta^{3}({\bf{x-y}})
\end{equation}
The formula that relates thecanonical conjugate of the density to 
 currents and particle densities may also be written down.
\begin{equation}
\Pi({\bf{x}}) = \int^{ {\bf{x}} }\mbox{ }d{\bf{l}}\mbox{ }
(-1/\rho({\bf{y}})){\bf{J}}({\bf{y}}) + \Phi([\rho];{\bf{x}})
 - \int^{ {\bf{x}} }\mbox{ }d{\bf{l}}\mbox{ }
[-i\mbox{ }\Phi, \nabla\Pi]({\bf{y}})
\end{equation}
The formula for the annhilation operator in terms of currents and densities is 
\begin{equation}
\psi({\bf{x}}) = e^{-i\int^{ {\bf{x}} }\mbox{ }d{\bf{l}}\mbox{ }
(-1/\rho({\bf{y}})){\bf{J}}({\bf{y}}) -i \Phi([\rho];{\bf{x}})
 +i \int^{ {\bf{x}} }\mbox{ }d{\bf{l}}\mbox{ }
[-i\mbox{ }\Phi, \nabla\Pi]({\bf{y}}) }e^{i\mbox{ }\Phi([\rho];{\bf{x}}) }
(\rho({\bf{x}}))^{\frac{1}{2}}
\end{equation} 
The solution for the phase functional $ \Phi $ has been given
 in an earlier preprint\cite{setlurexact}.
 In the high density limit, the answer may be
 written down as follows,
\begin{equation}
\Phi([\rho]; {\bf{x}}) = 0 
\end{equation}
in the bose case 
and in the fermi case,
\begin{equation}
\Phi([\rho]; {\bf{x}}) = \sum_{ {\bf{q}} \neq {\bf{0}} }
\rho_{ {\bf{q}} }U_{ {\bf{q}} }({\bf{x}})
\end{equation}
and,
\begin{equation}
U_{ {\bf{q}} }({\bf{x}}) = e^{-i \mbox{ }{\bf{q.x}} }U_{0}({\bf{q}})
\end{equation}
where again,
\begin{equation}
U_{0}({\bf{q}}) =  \frac{1}{N}(\frac{ \theta(k_{f} - |{\bf{q}}|)
- w_{1}({\bf{q}})  }{ w_{2}({\bf{q}}) })^{\frac{1}{2}}
\end{equation}
and,
\begin{equation}
w_{1}({\bf{q}}) = (\frac{1}{4\mbox{ }N\mbox{ }\epsilon_{ {\bf{q}} }^{2} })
\sum_{ {\bf{k}} }(\frac{ {\bf{k.q}} }{m})^{2}
(\Lambda_{ {\bf{k}} }({-{\bf{q}} }))^{2}
\end{equation}
\begin{equation}
w_{2}({\bf{q}}) = (\frac{1}{N})
\sum_{ {\bf{k}} }
(\Lambda_{ {\bf{k}} }({-{\bf{q}} }))^{2}
\end{equation}
here, 
$ \Lambda_{ {\bf{k}} }({\bf{q}}) = n_{F}({\bf{k}}+{\bf{q}}/2)(1- n_{F}({\bf{k}}-{\bf{q}}/2)) $

\end{document}